\documentclass{amsart}
\usepackage{amsfonts}

\usepackage{graphicx}
\usepackage{amscd}

\newtheorem{theorem}{Theorem}
\theoremstyle{plain}

\newtheorem{corollary}{Corollary}

\newtheorem{remark}{Remark}

\numberwithin{equation}{section}

\input{tcilatex}

\begin{document}
\title[Trapezoid Inequality]{A Generalised Trapezoid Type Inequality for Convex Functions}
\author{S.S. Dragomir}
\address{School of Communications and Informatics\\
Victoria University of Technology\\
PO Box 14428\\
Melbourne City MC\\
8001, Victoria, Australia.}
\email{sever@matilda.vu.edu.au}
\urladdr{http://rgmia.vu.edu.au/SSDragomirWeb.html}
\date{June 26, 2001.}

\begin{abstract}
A generalised trapezoid inequality for convex functions and applications for
quadrature rules are given. A refinement and a counterpart result for the
Hermite-Hadamard inequalities are obtained and some inequalities for pdf's
and $\left( HH\right) -$divergence measure are also mentioned.
\end{abstract}

\maketitle

\section{Introduction}

The following integral inequality for the generalised trapezoid formula was
obtained in \cite{2b} (see also \cite[p. 68]{1b}):

\begin{theorem}
\label{ta.1}Let $f:\left[ a,b\right] \rightarrow \mathbb{R}$ be a function
of bounded variation. We have the inequality 
\begin{eqnarray}
&&\left| \int_{a}^{b}f\left( t\right) dt-\left[ \left( x-a\right) f\left(
a\right) +\left( b-x\right) f\left( b\right) \right] \right|  \label{a.1} \\
&\leq &\left[ \frac{1}{2}\left( b-a\right) +\left| x-\frac{a+b}{2}\right| %
\right] \bigvee_{a}^{b}\left( f\right) ,  \notag
\end{eqnarray}
holding for all $x\in \left[ a,b\right] ,$ where $\bigvee_{a}^{b}\left(
f\right) $ denotes the total variation of $f$ on the interval $\left[ a,b%
\right] $.\newline
The constant $\frac{1}{2}$ is the best possible one.
\end{theorem}

This result may be improved if one assumes the monotonicity of $f$ as
follows (see \cite[p. 76]{1b})

\begin{theorem}
\label{ta.2}Let $f:\left[ a,b\right] \rightarrow \mathbb{R}$ be a monotonic
nondecreasing function on $\left[ a,b\right] $. Then we have the inequality: 
\begin{eqnarray}
&&\left| \int_{a}^{b}f\left( t\right) dt-\left[ \left( x-a\right) f\left(
a\right) +\left( b-x\right) f\left( b\right) \right] \right|  \label{a.2} \\
&\leq &\left( b-x\right) f\left( b\right) -\left( x-a\right) f\left(
a\right) +\int_{a}^{b}sgn\left( x-t\right) f\left( t\right) dt  \notag \\
&\leq &\left( x-a\right) \left[ f\left( x\right) -f\left( a\right) \right]
+\left( b-x\right) \left[ f\left( b\right) -f\left( x\right) \right]  \notag
\\
&\leq &\left[ \frac{1}{2}\left( b-a\right) +\left| x-\frac{a+b}{2}\right| %
\right] \left[ f\left( b\right) -f\left( a\right) \right]  \notag
\end{eqnarray}
for all $x\in \left[ a,b\right] $.\newline
The above inequalities are sharp.
\end{theorem}

If the mapping is Lipschitzian, then the following result holds as well \cite
{3b} (see also \cite[p. 82]{1b}).

\begin{theorem}
\label{ta.3}Let $f:\left[ a,b\right] \rightarrow \mathbb{R}$ be an $L-$%
Lipschitzian function on $\left[ a,b\right] ,$ i.e.., $f$ satisfies the
condition: 
\begin{equation}
\left| f\left( s\right) -f\left( t\right) \right| \leq L\left| s-t\right| 
\text{ \ for any \ }s,t\in \left[ a,b\right] \;\;\;\text{(}L>0\text{ is
given).}  \tag{L}  \label{L}
\end{equation}
Then we have the inequality: 
\begin{eqnarray}
&&\left| \int_{a}^{b}f\left( t\right) dt-\left[ \left( x-a\right) f\left(
a\right) +\left( b-x\right) f\left( b\right) \right] \right|  \label{a.3} \\
&\leq &\left[ \frac{1}{4}\left( b-a\right) ^{2}+\left( x-\frac{a+b}{2}%
\right) ^{2}\right] L  \notag
\end{eqnarray}
for any $x\in \left[ a,b\right] $.\newline
The constant $\frac{1}{4}$ is best in (\ref{a.3}).
\end{theorem}

If we would assume absolute continuity for the function $f$, then the
following estimates in terms of the Lebesgue norms of the derivative $%
f^{\prime }$ hold \cite[p. 93]{1b}.

\begin{theorem}
\label{ta.4}Let $f:\left[ a,b\right] \rightarrow \mathbb{R}$ be an
absolutely continuous function on $\left[ a,b\right] $. Then for any $x\in %
\left[ a,b\right] $, we have 
\begin{eqnarray}
&&\left| \int_{a}^{b}f\left( t\right) dt-\left[ \left( x-a\right) f\left(
a\right) +\left( b-x\right) f\left( b\right) \right] \right|  \label{a.4} \\
&\leq &\left\{ 
\begin{array}{lll}
\left[ \dfrac{1}{4}\left( b-a\right) ^{2}+\left( x-\dfrac{a+b}{2}\right) ^{2}%
\right] \left\| f^{\prime }\right\| _{\infty } & \text{if} & f^{\prime }\in
L_{\infty }\left[ a,b\right] ; \\ 
&  &  \\ 
\dfrac{1}{\left( q+1\right) ^{\frac{1}{q}}}\left[ \left( x-a\right)
^{q+1}+\left( b-x\right) ^{q+1}\right] ^{\frac{1}{q}}\left\| f^{\prime
}\right\| _{p} & \text{if} & f^{\prime }\in L_{p}\left[ a,b\right] , \\ 
&  & p>1,\;\frac{1}{p}+\frac{1}{q}=1; \\ 
\left[ \dfrac{1}{2}\left( b-a\right) +\left| x-\dfrac{a+b}{2}\right| \right]
\left\| f^{\prime }\right\| _{1}, &  & 
\end{array}
\right.  \notag
\end{eqnarray}
where $\left\| \cdot \right\| _{p}$ $\left( p\in \left[ 1,\infty \right]
\right) $ are the Lebesgue norms, i.e., 
\begin{equation*}
\left\| f^{\prime }\right\| _{\infty }=ess\sup\limits_{s\in \left[ a,b\right]
}\left| f^{\prime }\left( s\right) \right|
\end{equation*}
and 
\begin{equation*}
\left\| f^{\prime }\right\| _{p}:=\left( \int_{a}^{b}\left| f^{\prime
}\left( s\right) \right| ds\right) ^{\frac{1}{p}},\;\;p\geq 1.
\end{equation*}
\end{theorem}

In this paper we point out some similar results for convex functions.
Applications for quadrature formulae, for probability density functions and $%
HH-$Divergences in Information Theory are also considered.

\section{The Results}

The following theorem providing a lower bound for the difference 
\begin{equation*}
\left( x-a\right) f\left( a\right) +\left( b-x\right) f\left( b\right)
-\int_{a}^{b}f\left( t\right) dt
\end{equation*}
holds.

\begin{theorem}
\label{t1}Let $f:\left[ a,b\right] \rightarrow \mathbb{R}$ be a convex
function on $\left[ a,b\right] .$ Then for any $x\in \left( a,b\right) $ we
have the inequality 
\begin{eqnarray}
&&\frac{1}{2}\left[ \left( b-x\right) ^{2}f_{+}^{\prime }\left( x\right)
-\left( x-a\right) ^{2}f_{-}^{\prime }\left( x\right) \right]  \label{2.1} \\
&\leq &\left( x-a\right) f\left( a\right) +\left( b-x\right) f\left(
b\right) -\int_{a}^{b}f\left( t\right) dt.  \notag
\end{eqnarray}
The constant $\frac{1}{2}$ in the left hand side of (\ref{2.1}) is sharp in
the sense that it cannot be replaced by a larger constant.
\end{theorem}

\begin{proof}
It is easy to see that for any locally absolutely continuous function $%
f:\left( a,b\right) \rightarrow \mathbb{R}$, we have the identity 
\begin{equation}
\left( x-a\right) f\left( a\right) +\left( b-x\right) f\left( b\right)
-\int_{a}^{b}f\left( t\right) dt=\int_{a}^{b}\left( t-x\right) f^{\prime
}\left( t\right) dt  \label{2.2}
\end{equation}
for any $x\in \left( a,b\right) ,$ where $f^{\prime }$ is the derivative of $%
f$ which exists a.e. on $\left[ a,b\right] .$

Since $f$ is convex, then it is locally Lipschitzian and thus (\ref{2.2})
holds. Moreover, for any $x\in \left( a,b\right) ,$ we have the
inequalities: 
\begin{equation}
f^{\prime }\left( t\right) \leq f_{-}^{\prime }\left( x\right) \text{ \ for
\ a.e. \ }t\in \left[ a,x\right]   \label{2.3}
\end{equation}
and 
\begin{equation}
f^{\prime }\left( t\right) \geq f_{+}^{\prime }\left( x\right) \text{ \ for
\ a.e. \ }t\in \left[ x,b\right] .  \label{2.4}
\end{equation}
If we multiply (\ref{2.3}) by $x-t\geq 0$, $t\in \left[ a,x\right] $ and
integrate on $\left[ a,x\right] $, we get 
\begin{equation}
\int_{a}^{x}\left( x-t\right) f^{\prime }\left( t\right) dt\leq \frac{1}{2}%
\left( x-a\right) ^{2}f_{-}^{\prime }\left( x\right)   \label{2.5}
\end{equation}
and if we multiply (\ref{2.4}) by $t-x\geq 0$, $t\in \left[ x,b\right] $ and
integrate on $\left[ x,b\right] ,$ we also have 
\begin{equation}
\int_{x}^{b}\left( t-x\right) f^{\prime }\left( t\right) dt\geq \frac{1}{2}%
\left( b-x\right) ^{2}f_{+}^{\prime }\left( x\right) .  \label{2.6}
\end{equation}
Finally, if we subtract (\ref{2.5}) from (\ref{2.6}) and use the
representation (\ref{2.2}), we deduce the desired inequality (\ref{2.1}).

Now, assume that (\ref{2.1}) holds with a constant $C>0$ instead of $\frac{1%
}{2},$ i.e., 
\begin{eqnarray}
&&C\left[ \left( b-x\right) ^{2}f_{+}^{\prime }\left( x\right) -\left(
x-a\right) ^{2}f_{-}^{\prime }\left( x\right) \right]  \label{2.7} \\
&\leq &\left( x-a\right) f\left( a\right) +\left( b-x\right) f\left(
b\right) -\int_{a}^{b}f\left( t\right) dt.  \notag
\end{eqnarray}
Consider the convex function $f_{0}\left( t\right) :=k\left| t-\frac{a+b}{2}%
\right| ,$ $k>0$, $t\in \left[ a,b\right] .$ Then 
\begin{eqnarray*}
f_{0^{+}}^{\prime }\left( \frac{a+b}{2}\right) &=&k,\;\;\;f_{0^{-}}^{\prime
}\left( \frac{a+b}{2}\right) =-k, \\
f_{0}\left( a\right) &=&\frac{k\left( b-a\right) }{2}=f_{0}\left( b\right)
,\;\;\;\int_{a}^{b}f_{0}\left( t\right) dt=\frac{1}{4}k\left( b-a\right)
^{2}.
\end{eqnarray*}
If in (\ref{2.7}) we choose $f_{0}$ as above and $x=\frac{a+b}{2},$ then we
get 
\begin{equation*}
C\left[ \frac{1}{4}\left( b-a\right) ^{2}k+\frac{1}{4}\left( b-a\right) ^{2}k%
\right] \leq \frac{1}{4}k\left( b-a\right) ^{2}
\end{equation*}
giving $C\leq \frac{1}{2}$, and the sharpness of the constant is proved.
\end{proof}

Now, recall that the following inequality which is well known in the
literature as the \textit{Hermite-Hadamard inequality} for convex functions
holds 
\begin{equation}
f\left( \frac{a+b}{2}\right) \leq \frac{1}{b-a}\int_{a}^{b}f\left( t\right)
dt\leq \frac{f\left( a\right) +f\left( b\right) }{2}.  \tag{H-H}  \label{HH}
\end{equation}

The following corollary gives a sharp lower bound for the difference 
\begin{equation*}
\frac{f\left( a\right) +f\left( b\right) }{2}-\frac{1}{b-a}%
\int_{a}^{b}f\left( t\right) dt.
\end{equation*}

\begin{corollary}
\label{c1}Let $f:\left[ a,b\right] \rightarrow \mathbb{R}$ be a convex
function on $\left[ a,b\right] $. Then 
\begin{eqnarray}
0 &\leq &\frac{1}{8}\left[ f_{+}^{\prime }\left( \frac{a+b}{2}\right)
-f_{-}^{\prime }\left( \frac{a+b}{2}\right) \right] \left( b-a\right)
\label{2.8} \\
&\leq &\frac{f\left( a\right) +f\left( b\right) }{2}-\frac{1}{b-a}%
\int_{a}^{b}f\left( t\right) dt.  \notag
\end{eqnarray}
The constant $\frac{1}{8}$ is sharp.
\end{corollary}

The proof is obvious by the above theorem. The sharpness of the constant is
obtained for $f_{0}\left( t\right) =k\left| t-\frac{a+b}{2}\right| ,$ $t\in %
\left[ a,b\right] ,$ $k>0.$

When $x$ is a point of differentiability, we may state the following
corollary as well.

\begin{corollary}
\label{c2}Let $f$ be as in Theorem \ref{t1}. If $x\in \left( a,b\right) $ is
a point of differentiability for $f,$ then 
\begin{equation}
\left( b-a\right) \left( \frac{a+b}{2}-x\right) f^{\prime }\left( x\right)
\leq \left( x-a\right) f\left( a\right) +\left( b-x\right) f\left( b\right)
-\int_{a}^{b}f\left( t\right) dt.  \label{2.9}
\end{equation}
\end{corollary}

\begin{remark}
\label{r1}If $f:I\subseteq \mathbb{R\rightarrow R}$ is convex on $I$ and if
we choose $x\in $\r{I} (\r{I} is the interior of $I$), $b=x+\frac{h}{2},$ $%
a=x-\frac{h}{2},$ $h>0$ is such that $a,b\in I,$ then from (\ref{2.1}) we
may write 
\begin{equation}
0\leq \frac{1}{8}h^{2}\left[ f_{+}^{\prime }\left( x\right) -f_{-}^{\prime
}\left( x\right) \right] \leq \frac{f\left( a\right) +f\left( b\right) }{2}%
\cdot h-\int_{x-\frac{h}{2}}^{x+\frac{h}{2}}f\left( t\right) dt  \label{2.10}
\end{equation}
and the constant $\frac{1}{8}$ is sharp in (\ref{2.10}).
\end{remark}

The following result providing an upper bound for the difference 
\begin{equation*}
\left( x-a\right) f\left( a\right) +\left( b-x\right) f\left( b\right)
-\int_{a}^{b}f\left( t\right) dt
\end{equation*}
also holds.

\begin{theorem}
\label{t3}Let $f:\left[ a,b\right] \rightarrow \mathbb{R}$ \ be a convex
function on $\left[ a,b\right] $. Then for any $x\in \left[ a,b\right] ,$ we
have the inequality: 
\begin{eqnarray}
&&\left( x-a\right) f\left( a\right) +\left( b-x\right) f\left( b\right)
-\int_{a}^{b}f\left( t\right) dt  \label{2.11} \\
&\leq &\frac{1}{2}\left[ \left( b-x\right) ^{2}f_{-}^{\prime }\left(
b\right) -\left( x-a\right) ^{2}f_{+}^{\prime }\left( a\right) \right] . 
\notag
\end{eqnarray}
The constant $\frac{1}{2}$ is sharp in the sense that it cannot be replaced
by a smaller constant.
\end{theorem}

\begin{proof}
If either $f_{+}^{\prime }\left( a\right) =-\infty $ or $f_{-}^{\prime
}\left( b\right) =+\infty ,$ then the inequality (\ref{2.11}) evidently
holds true.

Assume that $f_{+}^{\prime }\left( a\right) $ and $f_{-}^{\prime }\left(
b\right) $ are finite.

Since $f$ is convex on $\left[ a,b\right] ,$ we have 
\begin{equation}
f^{\prime }\left( t\right) \geq f_{+}^{\prime }\left( a\right) \text{ \ for
\ a.e. \ }t\in \left[ a,x\right]   \label{2.12}
\end{equation}
and 
\begin{equation}
f^{\prime }\left( t\right) \leq f_{-}^{\prime }\left( b\right) \text{ \ for
\ a.e. \ }t\in \left[ x,b\right] .  \label{2.13}
\end{equation}
If we multiply (\ref{2.12}) by $\left( x-t\right) \geq 0,$ $t\in \left[ a,x%
\right] $ and integrate on $\left[ a,x\right] ,$ then we deduce 
\begin{equation}
\int_{a}^{x}\left( x-t\right) f^{\prime }\left( t\right) dt\geq \frac{1}{2}%
\left( x-a\right) ^{2}f_{+}^{\prime }\left( a\right)   \label{2.14}
\end{equation}
and if we multiply (\ref{2.13}) by $t-x\geq 0,$ $t\in \left[ x,b\right] $
and integrate on $\left[ x,b\right] ,$ then we also have 
\begin{equation}
\int_{x}^{b}\left( t-x\right) f^{\prime }\left( t\right) dt\leq \frac{1}{2}%
\left( b-x\right) ^{2}f_{-}^{\prime }\left( b\right) .  \label{2.15}
\end{equation}
Finally, if we subtract (\ref{2.14}) from (\ref{2.15}) and use the
representation (\ref{2.2}), we deduce the desired inequality (\ref{2.11}).

Now, assume that (\ref{2.11}) holds with a constant $D>0$ instead of $\frac{1%
}{2},$ i.e., 
\begin{eqnarray}
&&\left( x-a\right) f\left( a\right) +\left( b-x\right) f\left( b\right)
-\int_{a}^{b}f\left( t\right) dt  \label{2.16} \\
&\leq &D\left[ \left( b-x\right) ^{2}f_{-}^{\prime }\left( b\right) -\left(
x-a\right) ^{2}f_{+}^{\prime }\left( a\right) \right] .  \notag
\end{eqnarray}
If we consider the convex function $f_{0}:\left[ a,b\right] \rightarrow 
\mathbb{R}$, $f_{0}\left( t\right) =k\left| t-\frac{a+b}{2}\right| ,$ then
we have $f_{-}^{\prime }\left( b\right) =k,$ $f_{+}^{\prime }\left( a\right)
=-k$ and by (\ref{2.16}) we deduce for $x=\frac{a+b}{2}$ that 
\begin{equation*}
\frac{1}{4}k\left( b-a\right) ^{2}\leq D\left[ \frac{1}{4}k\left( b-a\right)
^{2}+\frac{1}{4}k\left( b-a\right) ^{2}\right]
\end{equation*}
giving $D\geq \frac{1}{2}$, and the sharpness of the constant is proved.
\end{proof}

The following corollary related to the Hermite-Hadamard inequality is
interesting as well.

\begin{corollary}
\label{c3}Let $f:\left[ a,b\right] \rightarrow \mathbb{R}$ be convex on $%
\left[ a,b\right] .$ Then 
\begin{equation}
0\leq \frac{f\left( a\right) +f\left( b\right) }{2}-\frac{1}{b-a}%
\int_{a}^{b}f\left( t\right) dt\leq \frac{1}{8}\left[ f_{-}^{\prime }\left(
b\right) -f_{+}^{\prime }\left( a\right) \right] \left( b-a\right) 
\label{2.17}
\end{equation}
and the constant $\frac{1}{8}$ is sharp.
\end{corollary}

\begin{remark}
\label{r2}Denote $B:=f_{-}^{\prime }\left( b\right) ,$ $A:=f_{+}^{\prime
}\left( a\right) $ and assume that $B\neq A,$ i.e., $f$ is not constant on $%
\left( a,b\right) .$ Then 
\begin{equation*}
\left( b-x\right) ^{2}B-\left( x-a\right) ^{2}A=\left( B-A\right) \left[
x-\left( \frac{bB-aA}{B-A}\right) \right] ^{2}-\frac{AB}{B-A}\left(
b-a\right) ^{2}
\end{equation*}
and by (\ref{2.11}) we get 
\begin{eqnarray}
&&\left( x-a\right) f\left( a\right) +\left( b-x\right) f\left( b\right)
-\int_{a}^{b}f\left( t\right) dt  \label{2.18} \\
&\leq &\left( B-A\right) \left[ x-\left( \frac{bB-aA}{B-A}\right) \right]
^{2}-\frac{AB}{\left( B-A\right) ^{2}}\left( b-a\right) ^{2}  \notag
\end{eqnarray}
for any $x\in \left[ a,b\right] .$\newline
If $A\geq 0,$ then $x_{0}=\frac{bB-aA}{B-A}\in \left[ a,b\right] $, and by (%
\ref{2.18}) for $x=\frac{bB-aA}{B-A}$ we get that 
\begin{equation}
0\leq \frac{1}{2}\cdot \frac{AB}{B-A}\left( b-a\right) \leq \frac{Bf\left(
a\right) -Af\left( b\right) }{B-A}-\frac{1}{b-a}\int_{a}^{b}f\left( t\right)
dt  \label{2.19}
\end{equation}
which is an interesting inequality in itself as well.
\end{remark}

\section{The Composite Case}

Consider the division $I_{n}:a=x_{0}<x_{1}<\dots <x_{n-1}<x_{n}=b$ and
denote $h_{i}:=x_{i+1}-x_{i}$ \ $\left( i=\overline{0,n-1}\right) .$ If $\xi
_{i}\in \left[ x_{i},x_{i+1}\right] $ $\left( i=\overline{0,n-1}\right) $
are intermediate points, then we will denote by 
\begin{equation}
G_{n}\left( f;I_{n},\mathbf{\xi }\right) :=\sum_{i=0}^{n-1}\left[ \left( \xi
_{i}-x_{i}\right) f\left( x_{i}\right) +\left( x_{i+1}-\xi _{i}\right)
f\left( x_{i+1}\right) \right]  \label{3.1}
\end{equation}
the generalised trapezoid rule associated to $f,$ $I_{n}$ and $\mathbf{\xi }%
. $

The following theorem providing upper and lower bounds for the remainder in
approximating the integral $\int_{a}^{b}f\left( t\right) dt$ of a convex
function $f$ in terms of the generalised trapezoid rule holds.

\begin{theorem}
\label{t4}Let $f:\left[ a,b\right] \rightarrow \mathbb{R}$ be a convex
function and $I_{n}$ and $\mathbf{\xi }$ be as above. Then we have: 
\begin{equation}
\int_{a}^{b}f\left( t\right) dt=G_{n}\left( f;I_{n},\mathbf{\xi }\right)
-S_{n}\left( f;I_{n},\mathbf{\xi }\right) ,  \label{3.2}
\end{equation}
where $G_{n}\left( f;I_{n},\mathbf{\xi }\right) $ is the generalised
Trapezoid Rule defined by (\ref{3.1}) and the remainder $S_{n}\left( f;I_{n},%
\mathbf{\xi }\right) $ satisfies the estimate: 
\begin{eqnarray}
&&\frac{1}{2}\left[ \sum_{i=0}^{n-1}\left( x_{i+1}-\xi _{i}\right)
^{2}f_{+}^{\prime }\left( \xi _{i}\right) -\sum_{i=0}^{n-1}\left( \xi
_{i}-x_{i}\right) ^{2}f_{-}^{\prime }\left( \xi _{i}\right) \right]
\label{3.3} \\
&\leq &S_{n}\left( f;I_{n},\mathbf{\xi }\right)  \notag \\
&\leq &\frac{1}{2}\left[ \left( b-\xi _{n-1}\right) ^{2}f_{-}^{\prime
}\left( b\right) +\sum_{i=1}^{n-1}\left[ \left( x_{i}-\xi _{i-1}\right)
^{2}f_{-}^{\prime }\left( x_{i}\right) -\left( \xi _{i}-x_{i}\right)
^{2}f_{+}^{\prime }\left( x_{i}\right) \right] \right.  \notag \\
&&-\left( \xi _{0}-a\right) ^{2}f_{+}^{\prime }\left( a\right) \bigg]. 
\notag
\end{eqnarray}
\end{theorem}

\begin{proof}
If we write the inequalities (\ref{2.1}) and (\ref{2.11}) on the interval $%
\left[ x_{i},x_{i+1}\right] $ and for the intermediate points $\xi _{i}\in %
\left[ x_{i},x_{i+1}\right] ,$ then we have 
\begin{eqnarray*}
&&\frac{1}{2}\left[ \left( x_{i+1}-\xi _{i}\right) ^{2}f_{+}^{\prime }\left(
x_{i}\right) -\left( \xi _{i}-x_{i}\right) ^{2}f_{-}^{\prime }\left( \xi
_{i}\right) \right] \\
&\leq &\left( \xi _{i}-x_{i}\right) f\left( x_{i}\right) +\left( x_{i+1}-\xi
_{i}\right) f\left( x_{i+1}\right) -\int_{x_{i}}^{x_{i+1}}f\left( t\right) dt
\\
&\leq &\frac{1}{2}\left[ \left( x_{i+1}-\xi _{i}\right) ^{2}f_{-}^{\prime
}\left( x_{i+1}\right) -\left( \xi _{i}-x_{i}\right) ^{2}f_{+}^{\prime
}\left( x_{i}\right) \right] .
\end{eqnarray*}
Summing the above inequalities over $i$ from $0$ to $n-1,$ we deduce 
\begin{eqnarray}
&&\frac{1}{2}\sum_{i=0}^{n-1}\left[ \left( x_{i+1}-\xi _{i}\right)
^{2}f_{+}^{\prime }\left( \xi _{i}\right) -\left( \xi _{i}-x_{i}\right)
^{2}f_{-}^{\prime }\left( \xi _{i}\right) \right]  \label{3.4} \\
&\leq &G_{n}\left( f;I_{n},\mathbf{\xi }\right) -\int_{a}^{b}f\left(
t\right) dt  \notag \\
&\leq &\frac{1}{2}\left[ \sum_{i=0}^{n-1}\left( x_{i+1}-\xi _{i}\right)
^{2}f_{-}^{\prime }\left( x_{i+1}\right) -\sum_{i=0}^{n-1}\left( \xi
_{i}-x_{i}\right) ^{2}f_{+}^{\prime }\left( x_{i}\right) \right] .  \notag
\end{eqnarray}
However, 
\begin{eqnarray*}
\sum_{i=0}^{n-1}\left( x_{i+1}-\xi _{i}\right) ^{2}f_{-}^{\prime }\left(
x_{i+1}\right) &=&\left( b-\xi _{n-1}\right) ^{2}f_{-}^{\prime }\left(
b\right) +\sum_{i=0}^{n-2}\left[ \left( x_{i+1}-\xi _{i}\right)
^{2}f_{-}^{\prime }\left( x_{i+1}\right) \right] \\
&=&\left( b-\xi _{n-1}\right) ^{2}f_{-}^{\prime }\left( b\right)
+\sum_{i=1}^{n-1}\left( x_{i}-\xi _{i-1}\right) ^{2}f_{-}^{\prime }\left(
x_{i}\right)
\end{eqnarray*}
and 
\begin{equation*}
\sum_{i=0}^{n-1}\left( \xi _{i}-x_{i}\right) ^{2}f_{+}^{\prime }\left(
x_{i}\right) =\sum_{i=1}^{n-1}\left( \xi _{i}-x_{i}\right) ^{2}f_{+}^{\prime
}\left( x_{i}\right) +\left( \xi _{0}-a\right) ^{2}f_{+}^{\prime }\left(
a\right)
\end{equation*}
and then, by (\ref{3.4}), we deduce the desired estimate (\ref{3.3}).
\end{proof}

The following corollary may be useful in practical applications.

\begin{corollary}
\label{c4}Let $f:\left[ a,b\right] \rightarrow \mathbb{R}$ be a
differentiable convex function on $\left[ a,b\right] .$ Then we have the
representation (\ref{3.2}) and $S_{n}\left( f;I_{n},\mathbf{\xi }\right) $
satisfies the estimate: 
\begin{eqnarray}
&&\sum_{i=0}^{n-1}\left( \frac{x_{i}+x_{i+1}}{2}-\xi _{i}\right)
h_{i}f^{\prime }\left( \xi _{i}\right)  \label{3.5} \\
&\leq &S_{n}\left( f;I_{n},\mathbf{\xi }\right)  \notag \\
&\leq &\frac{1}{2}\bigg[\left( b-\xi _{n-1}\right) ^{2}f_{-}^{\prime }\left(
b\right) -\left( \xi _{0}-a\right) ^{2}f_{+}^{\prime }\left( a\right)  \notag
\\
&&+\left. \sum_{i=1}^{n-1}\left[ \left( x_{i}-\frac{\xi _{i}+\xi _{i-1}}{2}%
\right) \left( \xi _{i}-\xi _{i-1}\right) f^{\prime }\left( x_{i}\right) %
\right] \right] .  \notag
\end{eqnarray}
\end{corollary}

We may also consider the trapezoid quadrature rule: 
\begin{equation}
T_{n}\left( f;I_{n}\right) :=\sum_{i=0}^{n-1}\frac{f\left( x_{i}\right)
+f\left( x_{i+1}\right) }{2}\cdot h_{i}.  \label{3.6}
\end{equation}
Using the above results, we may state the following corollary.

\begin{corollary}
\label{c5}Assume that $f:\left[ a,b\right] \rightarrow \mathbb{R}$ is a
convex function on $\left[ a,b\right] $ and $I_{n}$ is a division as above.
Then we have the representation 
\begin{equation}
\int_{a}^{b}f\left( t\right) dt=T_{n}\left( f;I_{n}\right) -Q_{n}\left(
f;I_{n}\right)  \label{3.7}
\end{equation}
where $T_{n}\left( f;I_{n}\right) $ is the mid-point quadrature formula
given in \ (\ref{3.6}) and the remainder $Q_{n}\left( f;I_{n}\right) $
satisfies the estimates 
\begin{eqnarray}
0 &\leq &\frac{1}{8}\sum_{i=0}^{n-1}\left[ f_{+}^{\prime }\left( \frac{%
x_{i}+x_{i+1}}{2}\right) -f_{-}^{\prime }\left( \frac{x_{i}+x_{i+1}}{2}%
\right) \right] h_{i}^{2}  \label{3.8} \\
&\leq &Q_{n}\left( f;I_{n}\right) \leq \frac{1}{8}\sum_{i=0}^{n-1}\left[
f_{+}^{\prime }\left( x_{i+1}\right) -f_{-}^{\prime }\left( x_{i}\right) %
\right] h_{i}^{2}.  \notag
\end{eqnarray}
The constant $\frac{1}{8}$ is sharp in both inequalities.
\end{corollary}

\section{Applications for P.D.F.s}

Let $X$ be a random variable with the \textit{probability density function} $%
f:\left[ a,b\right] \subset \mathbb{R\rightarrow }[0,\mathbb{\infty )}$ and
with \textit{cumulative distribution function }$F\left( x\right) =\Pr \left(
X\leq x\right) .$

The following theorem holds.

\begin{theorem}
\label{t4.1}If $f:\left[ a,b\right] \subset \mathbb{R\rightarrow R}_{+}$ is
monotonically increasing on $\left[ a,b\right] $, then we have the
inequality: 
\begin{eqnarray}
&&\frac{1}{2}\left[ \left( b-x\right) ^{2}f_{+}\left( x\right) -\left(
x-a\right) ^{2}f_{-}\left( x\right) \right] +x  \label{4.1} \\
&\leq &E\left( X\right)  \notag \\
&\leq &\frac{1}{2}\left[ \left( b-x\right) ^{2}f_{+}\left( b\right) -\left(
x-a\right) ^{2}f_{-}\left( a\right) \right] +x  \notag
\end{eqnarray}
for any $x\in \left( a,b\right) ,$ where $f_{\pm }\left( \alpha \right) $
represent respectively the right and left limits of $f$ in $\alpha $ and $%
E\left( X\right) $ is the expectation of $X.$\newline
The constant $\frac{1}{2}$ is sharp in both inequalities.\newline
The second inequality also holds for $x=a$ or $x=b.$
\end{theorem}

\begin{proof}
Follows by Theorem \ref{t1} and \ref{t3} applied for the convex cdf function 
$F\left( x\right) =\int_{a}^{x}f\left( t\right) dt,\;x\in \left[ a,b\right] $
and taking into account that 
\begin{equation*}
\int_{a}^{b}F\left( x\right) dx=b-E\left( X\right) .
\end{equation*}
\end{proof}

Finally, we may state the following corollary in estimating the expectation
of $X.$

\begin{corollary}
With the above assumptions, we have 
\begin{eqnarray}
&&\frac{1}{8}\left[ f_{+}\left( \frac{a+b}{2}\right) -f_{-}\left( \frac{a+b}{%
2}\right) \right] \left( b-a\right) ^{2}+\frac{a+b}{2}  \label{4.2} \\
&\leq &E\left( X\right) \leq \frac{1}{8}\left[ f_{+}\left( b\right)
-f_{-}\left( a\right) \right] \left( b-a\right) ^{2}+\frac{a+b}{2}.  \notag
\end{eqnarray}
\end{corollary}

\section{Applications for $HH-$Divergence}

Assume that a set $\chi $ and the $\sigma -$finite measure $\mu $ are given.
Consider the set of all probability densities on $\mu $ to be 
\begin{equation}
\Omega :=\left\{ p|p:\Omega \rightarrow \mathbb{R},\;p\left( x\right) \geq
0,\;\int_{\chi }p\left( x\right) d\mu \left( x\right) =1\right\} .
\label{6.1}
\end{equation}
Csisz\'{a}r's $f-$divergence is defined as follows \cite{4b} 
\begin{equation}
D_{f}\left( p,q\right) :=\int_{\chi }p\left( x\right) f\left[ \frac{q\left(
x\right) }{p\left( x\right) }\right] d\mu \left( x\right) ,\;p,q\in \Omega ,
\label{6.2}
\end{equation}
where $f$ is convex on $\left( 0,\infty \right) $. It is assumed that $%
f\left( u\right) $ is zero and strictly convex at $u=1.$ By appropriately
defining this convex function, various divergences are derived.

In \cite{5b}, Shioya and Da-te introduced the generalised Lin-Wong $f-$%
divergence $D_{f}\left( p,\frac{1}{2}p+\frac{1}{2}q\right) $ and the
Hermite-Hadamard $\left( HH\right) $ divergence 
\begin{equation}
D_{HH}^{f}\left( p,q\right) :=\int_{\chi }\frac{p^{2}\left( x\right) }{%
q\left( x\right) -p\left( x\right) }\left( \int_{1}^{\frac{q\left( x\right) 
}{p\left( x\right) }}f\left( t\right) dt\right) d\mu \left( x\right)
,\;p,q\in \Omega ,  \label{6.3}
\end{equation}
and, by the use of the Hermite-Hadamard inequality for convex functions,
proved the following basic inequality 
\begin{equation}
D_{f}\left( p,\frac{1}{2}p+\frac{1}{2}q\right) \leq D_{HH}^{f}\left(
p,q\right) \leq \frac{1}{2}D_{f}\left( p,q\right) ,  \label{6.4}
\end{equation}
provided that $f$ is convex and normalised, i.e., $f\left( 1\right) =0.$

The following result in estimating the difference 
\begin{equation*}
\frac{1}{2}D_{f}\left( p,q\right) -D_{HH}^{f}\left( p,q\right)
\end{equation*}
holds.

\begin{theorem}
\label{t6.1}Let $f:[0,\infty )\rightarrow \mathbb{R}$ be a normalised convex
function and $p,q\in \Omega .$ Then we have the inequality: 
\begin{eqnarray}
0 &\leq &\frac{1}{8}\left[ D_{f_{+}^{\prime }\cdot \left| \frac{\cdot +1}{2}%
\right| }\left( p,q\right) -D_{f_{-}^{\prime }\cdot \left| \frac{\cdot +1}{2}%
\right| }\left( p,q\right) \right]  \label{6.5} \\
&\leq &\frac{1}{2}D_{f}\left( p,q\right) -D_{HH}^{f}\left( p,q\right)  \notag
\\
&\leq &\frac{1}{8}D_{f_{-}^{\prime }\cdot \left( \cdot -1\right) }\left(
p,q\right) .  \notag
\end{eqnarray}
\end{theorem}

\begin{proof}
Using the double inequality 
\begin{eqnarray*}
0 &\leq &\frac{1}{8}\left[ f_{+}^{\prime }\left( \frac{a+b}{2}\right)
-f_{-}^{\prime }\left( \frac{a+b}{2}\right) \right] \left| b-a\right|  \\
&\leq &\frac{f\left( a\right) +f\left( b\right) }{2}-\frac{1}{b-a}%
\int_{a}^{b}f\left( t\right) dt \\
&\leq &\frac{1}{8}\left[ f_{-}\left( b\right) -f_{+}^{\prime }\left(
a\right) \right] \left( b-a\right) 
\end{eqnarray*}
for the choices $a=1$, $b=\frac{q\left( x\right) }{p\left( x\right) },$ $%
x\in \chi ,$ multiplying with $p\left( x\right) \geq 0$ and integrating over 
$x$ on $\chi $ we get 
\begin{eqnarray*}
0 &\leq &\frac{1}{8}\int_{\chi }\left[ f_{+}^{\prime }\left( \frac{p\left(
x\right) +q\left( x\right) }{2p\left( x\right) }\right) -f_{-}^{\prime
}\left( \frac{p\left( x\right) +q\left( x\right) }{2p\left( x\right) }%
\right) \right] \left| q\left( x\right) -p\left( x\right) \right| d\mu
\left( x\right)  \\
&\leq &\frac{1}{2}D_{f}\left( p,q\right) -D_{HH}^{f}\left( p,q\right)  \\
&\leq &\frac{1}{8}\int_{\chi }\left[ f_{-}^{\prime }\left( \frac{q\left(
x\right) }{p\left( x\right) }\right) -f_{+}^{\prime }\left( 1\right) \right]
\left( q\left( x\right) -p\left( x\right) \right) d\mu \left( x\right) ,
\end{eqnarray*}
which is clearly equivalent to (\ref{6.5}).
\end{proof}

\begin{corollary}
\label{c6.2}With the above assumptions and if $f$ is differentiable on $%
\left( 0,\infty \right) ,$ then 
\begin{equation}
0\leq \frac{1}{2}D_{f}\left( p,q\right) -D_{HH}^{f}\left( p,q\right) \leq 
\frac{1}{8}D_{f^{\prime }\cdot \left( \cdot -1\right) }\left( p,q\right) .
\label{6.6}
\end{equation}
\end{corollary}

\end{document}